\documentclass[11pt,a4paper]{amsart}

\usepackage[T1]{fontenc}
\usepackage{lmodern}
\usepackage{amsmath,amssymb,amsthm}
\usepackage[numbers,sort&compress]{natbib}
\usepackage[a4paper,margin=28mm]{geometry}
\usepackage{microtype}
\usepackage{needspace}
\usepackage{xurl}
\usepackage[hidelinks]{hyperref}
\input{glyphtounicode}
\hypersetup{
  pdftitle={The Universal Frame Characterization of C*-Algebras},
  pdfauthor={Mohammad B. Asadi; Zahra Hassanpour-Yakhdani},
  pdfsubject={Frames in Hilbert C*-modules},
  pdfkeywords={Hilbert C*-modules, frames, compact operators, positive definite kernels}
}
\theoremstyle{plain}
\newtheorem{theorem}{Theorem}[section]
\newtheorem{proposition}[theorem]{Proposition}
\newtheorem{lemma}[theorem]{Lemma}

\theoremstyle{definition}
\newtheorem{definition}[theorem]{Definition}
\newtheorem{example}[theorem]{Example}
\theoremstyle{remark}

\numberwithin{equation}{section}

\newcommand{\K}{\mathcal K}

\newcommand{\ip}[2]{\langle #1,#2\rangle}
\newcommand{\dens}{\operatorname{dens}}

\begin{document}

\title[The universal frame characterization]{The Universal Frame Characterization of $C^\ast$-Algebras}

\author[M. B. Asadi]{Mohammad B. Asadi}
\address{School of Mathematics, Statistics and Computer Science, College of Science,
University of Tehran, Tehran, 14155-6455, Iran}
\email{mb.asadi@ut.ac.ir}

\author[Z. Hassanpour-Yakhdani]{Zahra Hassanpour-Yakhdani}
\address{School of Mathematics, Statistics and Computer Science, College of Science,
University of Tehran, Tehran, 14155-6455, Iran}
\email{z.hasanpour@ut.ac.ir}

\thanks{Corresponding author: Mohammad B. Asadi.}

\begin{abstract}
We show that a $C^\ast$-algebra $A$ has the property that every Hilbert
$A$-module admits a frame if and only if $A$ is a $c_0$-direct sum of
algebras of compact operators. This settles the characterization question
of Amini, Asadi, Elliott, and Khosravi (2017), arising from the
frame-existence problem of Frank and Larson (2002).
\end{abstract}

\maketitle

\section{Introduction}

Frames for countably generated Hilbert $C^\ast$-modules are well understood: every such module admits a standard Parseval frame \cite[Proposition~2.3 and Theorem~2.5]{ArambaBakic}. For arbitrary modules, the corresponding existence problem depends on the coefficient algebra. Frank and Larson \cite[Problem~8.1]{FL} asked in 2002 whether every Hilbert module over a unital $C^\ast$-algebra has a normalized tight frame. Li \cite[Theorem~1.1]{Li} gave a negative answer and characterized the commutative unital algebras with universal frame existence as precisely the finite-dimensional ones. Amini, Asadi, Elliott, and Khosravi obtained the nonunital commutative counterpart and formulated the general characterization problem \cite[Theorem~1.4 and Questions~1.2, 1.5]{AAEK}.

The positive direction is known: by the results of Baki\'c and Gulja\v{s} \cite[Theorems~1 and~4]{BG}, every Hilbert module over a $C^\ast$-algebra of compact operators admits a standard Parseval frame. We prove the converse without separability or $\sigma$-unitality assumptions. More precisely, for every $C^\ast$-algebra $A$, we obtain
\[
 \text{every Hilbert $A$-module admits a frame}
 \quad\Longleftrightarrow\quad
 A\cong\bigoplus_{j\in J}^{c_0}\K(H_j).
\]
Here a frame is allowed to have an arbitrary index set, with ultraweak convergence of its frame sums. The characterization remains valid when the universal quantifier is restricted to full modules, or to modules with at most $\aleph_1$ generators. For an algebra not of compact-operator type, we construct a full frameless module and a frameless submodule with minimal generating cardinal $\aleph_1$. The construction also determines the generating cardinal and density character of the full module.

The obstruction is built directly in the category of Hilbert modules. A nonzero accumulation point of the spectrum of a positive element gives, through continuous functional calculus, a positive definite kernel indexed by $\aleph_1$ binary sequences. Suitable states distinguish two aspects of the resulting module: its distinguished elements are represented by only finitely many classes at each approximating state, but form an uncountable orthonormal family at a weak$^*$ cluster state. An upper Bessel bound, together with countable supports of module elements, is then incompatible with a lower frame bound. The underlying separation idea is related to Li's construction \cite{Li}; a state-based countability argument also occurs in \cite[proof of Theorem~2.8]{AAEK}.

Earlier examples of frameless modules were constructed in \cite{AH2,AFH}; see also \cite{Fufaev2025} for locally unital algebras. Related investigations of compactness and frames in Hilbert modules appear in \cite{FT2026,Fufaev2026}. The present argument uses neither representation-theoretic hypotheses on $A$ nor additional set-theoretic assumptions. As a consequence of the classification and \cite[Theorem~1.6]{AAEK}, it also addresses the question concerning maximal abelian subalgebras raised in \cite[Question~1.7]{AAEK}.

\section{Preliminaries and the spectral criterion}

Throughout, $A$ denotes a $C^\ast$-algebra. All Hilbert $A$-modules are right modules, with inner products linear in the second variable. A module $E$ is full if $\overline{\operatorname{span}}\ip EE=A$. A subset generates $E$ if its $A$-linear span is norm dense. We write $\operatorname{gen}(E)$ for the least cardinality of a generating set and $\dens(E)$ for the density character of $E$. The unitization $\widetilde A$ is taken to be $A$ itself when $A$ is unital. Spectra and functional calculi are computed in $\widetilde A$ unless otherwise indicated. We refer to \cite{Lance} for standard facts concerning Hilbert modules.

\begin{definition}\label{def:frame}
A family $(u_i)_{i\in I}$ in a Hilbert $A$-module $E$ is a frame (with ultraweak summation) if there exist constants $0<C\le D<\infty$ such that, for every $x\in E$, the net
\[
 S_F(x)=\sum_{i\in F}\ip{x}{u_i}\ip{u_i}{x},
 \qquad F\subset I\ \text{finite},
\]
converges ultraweakly in $A^{**}$ and its limit satisfies
\begin{equation}\label{eq:frame}
 C\ip xx\ \le\ \sum_{i\in I}\ip{x}{u_i}\ip{u_i}{x}\ \le\ D\ip xx
\end{equation}
in $A^{**}$. The frame is standard if the sums converge in norm in $A$, and Parseval if both bounds can be taken to be $1$. A family is a \emph{Bessel family} (with ultraweak summation) if there exists $D>0$ such that $S_F(x)\le D\ip xx$ for every finite $F\subseteq I$ and every $x\in E$. This bound implies ultraweak convergence of the increasing net of partial sums in $A^{**}$.
\end{definition}

Sums over arbitrary index sets are understood as limits of nets of finite partial sums. This is the convention of \cite[Section 2]{AAEK}; it includes standard frames. In particular, ruling out frames in this sense rules out standard frames. The empty family is a standard Parseval frame for the zero module.

For a state $\omega$ of $A$, also write $\omega$ for its normal extension to $A^{**}$, and consider the seminorms
\[
 \|a\|_\omega=\omega(a^*a)^{1/2},\qquad
 p_\omega(x)=\omega(\ip xx)^{1/2}.
\]
Both arise from positive semidefinite scalar inner products; in particular, they satisfy the triangle inequality and are bounded by their respective norms.

The following elementary inequalities will be used repeatedly.
\begin{lemma}\label{lem:state}
If $(u_i)_{i\in I}$ is a Bessel family with upper bound $D$, then for every state $\omega$ and $x\in E$,
\begin{equation}\label{eq:stateupper}
 \sum_{i\in I}\|\ip{u_i}{x}\|_\omega^2\le Dp_\omega(x)^2.
\end{equation}
For fixed $\omega$ and $x$, at most countably many summands are nonzero. If $x,y\in E$ satisfy $p_\omega(x-y)=0$, then $\|\ip{u_i}{x-y}\|_\omega=0$ for every $i$.
If the family is a frame with lower bound $C$, then
\begin{equation}\label{eq:state}
 Cp_\omega(x)^2\le \sum_{i\in I}\|\ip{u_i}{x}\|_\omega^2.
\end{equation}
\end{lemma}
\begin{proof}
Apply the normal extension of $\omega$ to the bounded increasing net of partial Bessel sums; compare \cite[Proposition~3.1]{Li} in the unital setting. A summable family of nonnegative numbers has countable support. Applying the upper inequality to $x-y$ proves the vanishing assertion. The lower estimate for frames follows by applying $\omega$ to the lower frame inequality.
\end{proof}

We say that a $C^\ast$-algebra $A$ is \emph{of compact-operator type} if it satisfies either of the equivalent conditions in Proposition~\ref{prop:spectral}(a)--(b). We use the standard structural description of $C^\ast$-algebras of compact operators \cite[p.~251]{BG} and the spectral criterion for such algebras \cite[Theorem~5.3]{Roch}, applied to $A$ as an ideal in its unitization. The self-adjoint and positive-element formulations are equivalent by spectral mapping under the absolute-value function.

\begin{proposition}\label{prop:spectral}
The following conditions are equivalent:
\begin{enumerate}
\item[(a)] $A$ is isomorphic to a $C^\ast$-subalgebra of compact operators on a Hilbert space.
\item[(b)] $A\cong\bigoplus_{j\in J}^{c_0}\K(H_j)$ for a family of Hilbert spaces $(H_j)_{j\in J}$.
\item[(c)] For every $a\in A_+$, $\sigma(a)$ has no nonzero accumulation point.
\end{enumerate}
\end{proposition}

The next lemma supplies the required spectral states.
\begin{lemma}\label{lem:spectralstates}
If $A$ is not of compact-operator type, there exist positive contractions $b,e\in A$, numbers $r_n\in(0,1]$ tending to zero, states $\varphi_n$ of $A$, and a weak$^*$ cluster state $\varphi$ of $(\varphi_n)$ such that
\begin{equation}\label{eq:spectralstates}
 \varphi_n(e g(b)e)=g(r_n),\qquad
 \varphi(e g(b)e)=g(0)
\end{equation}
for every $g\in C([0,1])$. In particular, $\varphi_n(e^2)=\varphi(e^2)=1$.
\end{lemma}
\begin{proof}
By Proposition~\ref{prop:spectral}, choose $a\in A_+$ and distinct $\lambda_n\in\sigma(a)$ tending to $\lambda>0$, with $\lambda_n\ne\lambda$ and $\lambda_n>\lambda/2$. Put $f(t)=t(t-\lambda)^2$ on $\sigma(a)\cup\{0\}$ and
\[
 b=f(a)/\|f(a)\|,\qquad r_n=f(\lambda_n)/\|f(a)\|.
\]
Since $f(a)\ne0$, the denominator is nonzero. Moreover, $r_n>0$ and $r_n\to0$. As $f(0)=0$, we also have $b\in A$. Let $e=h(a)$, where $0\le h\le1$ is continuous, $h(0)=0$, and $h(t)=1$ for $t\ge\lambda/2$. Extend evaluation at $\lambda_n$ on $C^\ast(1,a)\subseteq\widetilde A$ to a state $\widetilde\varphi_n$ of $\widetilde A$. By weak$^*$ compactness of the state space of $\widetilde A$, a subnet converges to a state $\widetilde\varphi$. Each of these states, as well as their limit, takes the value $1$ on $e^2$. Consequently $\varphi_n=\widetilde\varphi_n|_A$ and $\varphi=\widetilde\varphi|_A$ are states, and the corresponding subnet converges weak$^*$ to $\varphi$ in $A^*$. Since $\widetilde\varphi_n$ agrees with evaluation at $\lambda_n$ on $C^\ast(1,a)$, continuous functional calculus gives
\[
\begin{aligned}
\varphi_n(e g(b)e)
 &= h(\lambda_n)^2
    g\!\left(\frac{f(\lambda_n)}{\|f(a)\|}\right)\\
 &=g(r_n)
\end{aligned}
\]
for every $g\in C([0,1])$. The second equality in \eqref{eq:spectralstates} follows by weak$^*$ convergence along the chosen subnet and continuity of $g$ at zero. The elements $eg(b)e$ belong to $A$, even when $g(b)$ is computed in the unitization.
\end{proof}

\section{The universal frame characterization}

Fix a subset $S\subseteq\{0,1\}^{\mathbb N}$ of cardinality $\aleph_1$ containing the countable set $S_0$ of all eventually zero sequences. Such a set exists in ZFC. For distinct $s,t\in S$, let $m(s,t)$ denote the length of their common initial segment, which may be zero. Put
\[
 q_n(r)=\min\{1,nr\}\quad(n\ge0,\ 0\le r\le1),\qquad
 g_{st}(r)=
 \begin{cases}
 1,&s=t,\\
 q_{m(s,t)}(r),&s\ne t .
 \end{cases}
\]
In particular, $g_{st}(0)=\delta_{st}$.

\begin{lemma}\label{lem:kernel}
The matrix $(g_{st})_{s,t\in F}$ belongs to $M_{|F|}(C([0,1]))_+$ for every nonempty finite subset $F$ of $S$.
\end{lemma}
\begin{proof}
For $n\ge1$, let $P_n$ be the matrix indexed by $F$ whose $(s,t)$-entry is $1$ when $s$ and $t$ share their first $n$ coordinates, and $0$ otherwise. Each $P_n$ is positive, being block diagonal with all-ones blocks after a permutation of the indices. Choose $N\ge1$ larger than every $m(s,t)$ for distinct $s,t\in F$. An entrywise calculation gives
\[
 (g_{st}(r))_{s,t\in F}
  =(1-q_N(r))I+
    \sum_{n=1}^N(q_n(r)-q_{n-1}(r))P_n .
\]
Indeed, the diagonal entries on the right-hand side equal $1-q_N(r)+\sum_{n=1}^N(q_n(r)-q_{n-1}(r))=1$. For distinct $s,t\in F$, the $(s,t)$-entry is $\sum_{n=1}^{m(s,t)}(q_n(r)-q_{n-1}(r))=q_{m(s,t)}(r)=g_{st}(r)$, with the empty sum understood as zero. All coefficients are nonnegative. Pointwise positivity of this continuous matrix-valued function proves the claim.
\end{proof}

The kernel in Lemma~\ref{lem:kernel} determines a Hilbert module with the required generating properties.
\begin{proposition}\label{prop:module}
For any positive contraction $b\in A$, there exists a full Hilbert $A$-module $E_b$ and adjointable isometries $v_s:A_A\to E_b$ ($s\in S$) satisfying
\begin{equation}\label{eq:gram}
 v_s(a)c=v_s(ac),\qquad
 \ip{v_s(a)}{v_t(c)}=a^*g_{st}(b)c.
\end{equation}
The linear span of $\bigcup_{s\in S}v_s(A)$ is dense in $E_b$. Every element of $E_b$ belongs to the closed linear span of $v_s(A)$ for $s$ in some countable subset of $S$. Furthermore,
\[
 \operatorname{gen}(E_b)\le\max\{\aleph_1,\operatorname{gen}(A_A)\},
 \qquad
 \dens(E_b)\le\max\{\aleph_1,\dens(A)\}.
\]
\end{proposition}
\begin{proof}
By Lemma~\ref{lem:kernel} and continuous functional calculus, $k(s,t)=L_{g_{st}(b)}$ defines a positive definite kernel with values in $\mathcal L(A_A)$, where $L_{g_{st}(b)}$ denotes adjointable left multiplication on $A_A$. Murphy's minimal Kolmogorov decomposition theorem \cite[Theorem~2.3]{Murphy} provides a Hilbert $A$-module $E_b$ and adjointable maps $V(s):A_A\to E_b$ such that $V(s)^*V(t)=k(s,t)$ and the linear span of $V(s)A$ over $s\in S$ is dense. Taking $v_s(a)=V(s)a$ gives \eqref{eq:gram}; the diagonal identity implies that each $v_s$ is isometric. Moreover, for every $d\in A_+$ and any fixed $s\in S$, we have
\[
 d=\ip{v_s(d^{1/2})}{v_s(d^{1/2})}.
\]
Since $A=\operatorname{span}(A_+)$, it follows that $\operatorname{span}\ip{E_b}{E_b}=A$. In particular, $E_b$ is full.

Given $u\in E_b$, by the density of the linear span of the $v_s(A)$, for each $n\geq1$ there exist a finite subset $F_n\subseteq S$ and elements $a_{n,s}\in A$ such that
\[
 \left\|u-\sum_{s\in F_n}v_s(a_{n,s})\right\|<\frac1n.
\]
The set $R=\bigcup_{n\geq1}F_n$ is at most countable, and hence
\[
 u\in\overline{\operatorname{span}}\{v_s(a):s\in R,\ a\in A\}.
\]
For the density bound, choose a norm-dense subset $D\subseteq A$ with $|D|=\dens(A)$. The finite $(\mathbb Q+i\mathbb Q)$-linear combinations of $\{v_s(d):s\in S,\ d\in D\}$ form a dense subset of $E_b$ of cardinality at most $\max\{\aleph_1,\dens(A)\}$. For the generating bound, let $T\subseteq A$ generate $A_A$. Then $\{v_s(t):s\in S,\ t\in T\}$ generates $E_b$, since each $v_s$ is right $A$-linear. Taking $|T|=\operatorname{gen}(A_A)$ gives the asserted bound.
\end{proof}

For $R\subseteq S$, write $E_b(R)=\overline{\operatorname{span}}\{v_t(a):t\in R,\ a\in A\}$, a closed submodule. We say that $R$ supports $u$ if $u\in E_b(R)$.

We now construct the modules that obstruct universal frame existence.
\begin{theorem}\label{thm:witness}
Let $A$ be a $C^\ast$-algebra not of compact-operator type. There exist a full Hilbert $A$-module $E_b$ and a closed submodule $E_0\subseteq E_b$ with the following properties:
\begin{enumerate}
\item Neither $E_b$ nor $E_0$ admits a frame, even with ultraweak summation.
\item The minimal generating cardinal of $E_0$ is $\aleph_1$.
\item The minimal generating cardinal and density character satisfy
\[
\operatorname{gen}(E_b)=\max\{\aleph_1,\operatorname{gen}(A_A)\},
\qquad \dens(E_b)=\max\{\aleph_1,\dens(A)\}.
\]
\end{enumerate}
\end{theorem}
\begin{proof}
Choose $b,e,r_n,\varphi_n,\varphi$ as in Lemma~\ref{lem:spectralstates}, and form $E_b$ by Proposition~\ref{prop:module}. For $s\in S$, put
\[
 x_s=v_s(e),\qquad
 E_0=\overline{\operatorname{span}}\{x_s a:s\in S,\ a\in A\}.
\]
For any approximate identity $(e_\lambda)$ of $A$, $x_se_\lambda\to x_s$; hence $x_s\in E_0$, and
\begin{equation}\label{eq:xgram}
 \ip{x_s}{x_t}=e g_{st}(b)e.
\end{equation}
Thus $p_{\varphi_n}(x_s)=p_\varphi(x_s)=1$ and $\varphi(\ip{x_s}{x_t})=\delta_{st}$.

To prove (1), let $(u_i)_{i\in I}$ be an arbitrary Bessel family in either module, with upper bound $D$. Set
\[
 I_0=\bigcup_{n\ge1,\ t\in S_0}
       \{i\in I:\|\ip{u_i}{x_t}\|_{\varphi_n}>0\}.
\]
Lemma~\ref{lem:state} shows that $I_0$ is countable. For each $n$, choose $N_n$ with $N_n r_n\ge1$. Given $s\in S$, choose $t\in S_0$ with the same first $N_n$ coordinates. Then $g_{st}(r_n)=1$ (also when $s=t$), and hence
\[
 p_{\varphi_n}(x_s-x_t)^2=2-2g_{st}(r_n)=0.
\]
For $i\notin I_0$, the Bessel upper bound applied to $x_s-x_t$ and the definition of $I_0$ give
\[
 \|\ip{u_i}{x_s}\|_{\varphi_n}
 \le \|\ip{u_i}{x_s-x_t}\|_{\varphi_n}
       +\|\ip{u_i}{x_t}\|_{\varphi_n}=0.
\]
Thus
\[
 \|\ip{u_i}{x_s}\|_{\varphi_n}=0
 \quad(n\ge1,\ s\in S,\ i\notin I_0).
\]
For each fixed $s,i$, the positive element $\ip{x_s}{u_i}\ip{u_i}{x_s}$ is in $A$. Evaluating this fixed element along the weak$^*$ convergent subnet of $(\varphi_n)$ gives
\begin{equation}\label{eq:outside}
 \|\ip{u_i}{x_s}\|_\varphi=0
 \quad(s\in S,\ i\notin I_0).
\end{equation}
Only fixed coefficients are passed to the limit; no interchange with an infinite sum is involved.

For each $i\in I_0$, choose a countable support $R_i\subseteq S$ for $u_i$ in $E_b$. Set $R=\bigcup_{i\in I_0}R_i$, a countable set. If $s\in S\setminus R$, $t\in R$ and $a\in A$, write $q=q_{m(t,s)}$. From \eqref{eq:gram} and the positivity of $aa^*$,
\begin{align*}
 \|\ip{v_t(a)}{x_s}\|_\varphi^2
 &=\varphi(e q(b)aa^*q(b)e)\\
 &\le\|a\|^2\varphi(e q(b)^2e)=\|a\|^2q(0)^2=0.
\end{align*}
The estimate $\|\ip{u}{x_s}\|_\varphi\le\|u\|\,\|x_s\|$ extends this vanishing to finite sums and their norm closure, hence to every $u\in E_b(R)$ and, in particular, to $u_i$ for $i\in I_0$. Together with \eqref{eq:outside}, this shows that for every $s\notin R$ and every $i\in I$,
\[
 \varphi\big(\ip{x_s}{u_i}\ip{u_i}{x_s}\big)=0.
\]

If the Bessel family were a frame with lower bound $C>0$, choose any $s\in S\setminus R$. The lower inequality in Lemma~\ref{lem:state} would give
\[
 0<C=Cp_\varphi(x_s)^2
 \le\sum_{i\in I}\|\ip{u_i}{x_s}\|_\varphi^2=0,
\]
a contradiction. This proves (1). Fullness of $E_b$ is part of Proposition~\ref{prop:module}.

The family $(x_s)_{s\in S}$ generates $E_0$, so $\operatorname{gen}(E_0)\le\aleph_1$. Neither module is countably generated, because neither admits a frame; every countably generated Hilbert module has a standard Parseval frame \cite[Proposition~2.3 and Theorem~2.5]{ArambaBakic}. Thus $\operatorname{gen}(E_0)=\aleph_1$ and $\operatorname{gen}(E_b)\ge\aleph_1$.

Proposition~\ref{prop:module} gives $\operatorname{gen}(E_b)\le\max\{\aleph_1,\operatorname{gen}(A_A)\}$. For any fixed $s$, the adjointable map $v_s^*:E_b\to A_A$ is surjective, since $v_s^*v_s=1_{A_A}$. Thus the image under $v_s^*$ of any generating set for $E_b$ generates $A_A$, giving $\operatorname{gen}(E_b)\ge\operatorname{gen}(A_A)$. Together with $\operatorname{gen}(E_b)\ge\aleph_1$, this proves the generating-cardinal equality in (3).

Finally, the $x_s$ are mutually separated by distance at least $\sqrt 2$, as witnessed by $p_\varphi$, and the maps $v_s$ embed $A$ isometrically. Hence $\dens(E_b)\ge\max\{\aleph_1,\dens(A)\}$. The reverse inequality follows from Proposition~\ref{prop:module}.
\end{proof}

We recall the following consequence of the orthonormal-basis theorem and Parseval's identity of Baki\'c and Gulja\v{s} \cite[Theorems~1 and~4]{BG}.

\begin{proposition}\label{prop:compact}
If $A$ is of compact-operator type, every Hilbert $A$-module has a standard Parseval frame.
\end{proposition}

\Needspace{19\baselineskip}
Combining this obstruction with the compact-operator case yields the characterization.
\begin{theorem}\label{thm:classification}
For an arbitrary $C^\ast$-algebra $A$, the following conditions are equivalent:
\begin{enumerate}
\item[(i)] $A$ is isomorphic to a $C^\ast$-subalgebra of compact operators; equivalently,
\[
 A\cong\bigoplus_{j\in J}^{c_0}\K(H_j).
\]
\item[(ii)] Every Hilbert $A$-module admits a standard Parseval frame.
\item[(iii)] Every Hilbert $A$-module admits a standard frame.
\item[(iv)] Every Hilbert $A$-module admits a frame with ultraweak summation.
\item[(v)] Every Hilbert $A$-module generated by at most $\aleph_1$ elements admits a frame with ultraweak summation.
\item[(vi)] Every full Hilbert $A$-module admits a frame with ultraweak summation.
\end{enumerate}
If $A$ is unital, these conditions are also equivalent to $A$ being finite-dimensional. If $A$ admits at least one full Hilbert $A$-module with at most $\aleph_1$ generators, the conditions are also equivalent to requiring a frame for every full Hilbert $A$-module with at most $\aleph_1$ generators. The latter existence hypothesis holds, in particular, if $A$ is $\sigma$-unital.
\end{theorem}
\begin{proof}
Proposition~\ref{prop:compact} gives (i)$\Rightarrow$(ii), while (ii)$\Rightarrow$(iii)$\Rightarrow$(iv)$\Rightarrow$(v) and (iv)$\Rightarrow$(vi) are immediate. If (i) fails, Theorem~\ref{thm:witness} gives a frameless module $E_0$ with $\operatorname{gen}(E_0)=\aleph_1$, contradicting (v), and a full frameless module $E_b$, contradicting (vi). Thus both (v) and (vi) imply (i). The structural equivalence in (i) is Proposition~\ref{prop:spectral}. A unital compact-operator algebra is finite-dimensional \cite[p.~2]{AAEK}.

For the additional equivalence, let $F$ be a full Hilbert $A$-module generated by at most $\aleph_1$ elements. If (i) fails, $E_0\oplus F$ is full and has at most $\aleph_1$ generators. A frame $((u_i,w_i))_{i\in I}$ on it would give a frame $(u_i)_{i\in I}$ on $E_0$: for $x\in E_0$, one has $\ip{(u_i,w_i)}{(x,0)}=\ip{u_i}{x}$, so the frame inequalities and convergence follow by testing on $(x,0)$. This contradicts Theorem~\ref{thm:witness}. Conversely, (ii) supplies a frame for every such module. If $A$ is $\sigma$-unital, $A_A$ is a full countably generated Hilbert module. The zero algebra is covered by the empty sum and the empty frame.
\end{proof}

\begin{example}
Every infinite-dimensional von Neumann algebra (viewed as a $C^\ast$-algebra), including $B(H)$ for infinite-dimensional $H$ and $\ell^\infty(J)$ for infinite $J$, fails the universal frame property. The same holds for infinite-dimensional UHF algebras, irrational rotation algebras, the Jiang--Su algebra, the Toeplitz algebra, and the Cuntz algebras $\mathcal O_n$ ($2\le n\le\infty$). On the other hand, $\K(H)$ and $c_0(J)$ have the property for arbitrary cardinalities. The property also fails for infinite-dimensional unitizations of compact-operator algebras; compare \cite{AH2,AFH}. For a locally compact Hausdorff space $X$, the theorem shows that every Hilbert $C_0(X)$-module admits a frame precisely when $X$ is discrete.
\end{example}

\section*{Data availability}
No datasets were generated or analyzed in this theoretical study.

\section*{Funding}
This research did not receive any specific grant from funding agencies in the public, commercial, or not-for-profit sectors.

\Needspace{9\baselineskip}
\section*{Declaration of generative AI and AI-assisted technologies in the manuscript preparation process}
During the preparation of this work, the authors used ChatGPT (OpenAI) for language editing, literature searches, mathematical discussions, and refinement of proofs. The mathematical ideas and arguments were developed through the authors' own contributions and critical discussions with the AI tool. After using this tool, the authors reviewed and edited the content as needed, independently examined and verified all proofs and mathematical arguments, and take full responsibility for the content of the published article.

\end{document}